# INTEGRATED FUNCTIONALS OF NORMAL AND FRACTIONAL PROCESSES[1]

By Boris Buchmann and Ngai Hang Chan

*Monash University and Chinese University of Hong Kong*

Consider $Z_t^f(u) = \int_0^{tu} f(N_s)\,ds$, $t > 0$, $u \in [0,1]$, where $N = (N_t)_{t \in \mathbb{R}}$ is a normal process and $f$ is a measurable real-valued function satisfying $Ef(N_0)^2 < \infty$ and $Ef(N_0) = 0$. If the dependence is sufficiently weak Hariz [*J. Multivariate Anal.* **80** (2002) 191–216] showed that $Z_t^f/t^{1/2}$ converges in distribution to a multiple of standard Brownian motion as $t \to \infty$. If the dependence is sufficiently strong, then $Z_t/(EZ_t(1)^2)^{1/2}$ converges in distribution to a higher order Hermite process as $t \to \infty$ by a result by Taqqu [*Wahrsch. Verw. Gebiete* **50** (1979) 53–83]. When passing from weak to strong dependence, a unique situation encompassed by neither results is encountered. In this paper, we investigate this situation in detail and show that the limiting process is still a Brownian motion, but a nonstandard norming is required. We apply our result to some functionals of fractional Brownian motion which arise in time series. For all Hurst indices $H \in (0,1)$, we give their limiting distributions. In this context, we show that the known results are only applicable to $H < 3/4$ and $H > 3/4$, respectively, whereas our result covers $H = 3/4$.

**1. Introduction.** With the increase popularity of using self-similar processes to model econometric time series (see Ballie [2]) and network traffic (see Willinger, Taqqu, Sherman and Wilson [27]), understanding the limiting behavior of integrated functionals of normal processes becomes a challenging and important task. Functional limit theorems of this nature have found applications in various disciplines (see, e.g., Bhattachary, Gupta and

Received April 2007; revised March 2008.

[1]Supported in part by HKSAR-RGC Grants CUHK400305 and CUHK400408 and ARC Grants DP0664603 and CE0348217.

Waymire [4] for an interesting application in hydrology). One of the early results of this nature is given in the seminal paper by Taqqu [25], where it was shown that the limit of such processes are in the domain of attraction of higher order Hermite processes when the underlying normal process exhibits strong dependence. This problem was subsequently studied in Hariz [15] when the underlying normal process is weakly dependent. The problem becomes more intriguing when transitions from weak to strong dependence take place. Among others, we establish the limiting distribution of integrated functionals of fractional Brownian motions for all Hurst exponents in $(0,1)$ in this paper.

Let $(\Omega, \mathcal{F}, P)$ be a probability space carrying a standard Brownian motion $W = (W_t)_{t \in \mathbb{R}}$ and let $a$ be a Borel function with $\int a^2(s)\,ds = 1$. We consider a normal process $N$ defined by

$$(1.1) \qquad N_t = \int_{\mathbb{R}} a(t-s)\,dW_s, \qquad t \in \mathbb{R},$$

and also its covariance function $r : \mathbb{R} \to \mathbb{R}$, where we set

$$r(t) = E(N_t N_0), \qquad t \geq 0.$$

For $k \in \mathbb{N}_0$, let $\mathcal{H}_k$ denote the $k$th Hermite polynomial with leading coefficient one defined by

$$(1.2) \qquad \mathcal{H}_k(x) = (-1)^k \exp\left(\frac{x^2}{2}\right) \frac{d^k}{dx^k} \exp\left(-\frac{x^2}{2}\right).$$

In the sequel, let $f$ be a measurable real-valued function such that $f$ satisfies $Ef(N_0) = 0$ and $Ef^2(N_0) < \infty$. Then $f$ can be expanded into an orthogonal series as

$$(1.3) \qquad f(x) = \sum_{k=1}^{\infty} \frac{c_k}{k!} \mathcal{H}_k(x),$$

where the coefficients are given by

$$(1.4) \qquad c_k = E\mathcal{H}_k(N_0) f(N_0), \qquad k \in \mathbb{N}.$$

The series in (1.3) converges in $L^2(\nu)$ where $\nu$ denotes the standard normal law on the Borel sets of $\mathbb{R}$. In this case, the constant

$$q = q(f) = \inf\{k : c_k \neq 0\}$$

is called the *Hermite rank* of $f$ by Taqqu in his pioneering paper [24].

The Hermite rank plays a crucial role in the asymptotic behavior of sample means of nonlinear functionals for both Gaussian processes and fields indexed by continuous and discrete parameters (cf. Breuer and Major [7], Chambers and Slud [11], Dobrushin and Major [13], Hariz [15] and Taqqu [24] and [25] and references therein).

In particular, Hariz [15] showed that, if

$$\int_0^\infty |r(s)|^q \, ds < \infty, \tag{1.5}$$

then a central limit theorem (CLT) holds, that is, the finite dimensional distributions of $Z_t^f(\cdot)/t^{1/2} = \int_0^{t\cdot} f(N_u) \, du/t^{1/2}$ converges to $\sigma W$, as $t \to \infty$, where

$$\sigma^2 = \sum_{k=q}^\infty \frac{c_k^2}{k!} \int r^k(u) \, du. \tag{1.6}$$

The integral in (1.5) is finite if the covariance function of $N$ tends to zero sufficiently fast, that is, the dependence of the underlying process $N$ is sufficiently weak. When the underlying process is strongly dependent, a noncentral limit theorem (NCLT) holds: Under some technical conditions on the kernel $a$, $Z^f$ is in the domain of attraction of a Hermite process of order $q$. This was shown by Taqqu [25]. However, there are situations where the quantity in (1.5) is infinite, but the dependence is still too weak to satisfy a NCLT. This situation typically occurs when one wants to explore the domain of attraction of integrated functionals of fractional Ornstein–Uhlenbeck processes for all Hurst indices $H \in (0,1)$. As a standard rule, the CLT covers $H \in (0, 1 - 1/(2q))$ whereas the NCLT is applicable to $H \in (1 - 1/(2q), 1)$ $[q \geq 2]$. In this paper, we are interested in the boundary corresponding to $H = 1 - 1/(2q)$.

The paper is organized as follows. In Section 2, we state a general theorem dealing with the boundary case. Also, we give an application of our result to fractional Ornstein–Uhlenbeck processes. An important application of results of this nature arises from the study of strongly dependent times series (cf. Buchmann and Chan [8]). In Section 3, we review and extend the asymptotic theory of the so-called unit-root problem as an application of Section 2. In this context, the CLT and the NCLT are only applicable to $H \in (0, 3/4)$ and $H \in (3/4, 1)$, respectively, whereas our theorem covers $H = 3/4$. All proofs are found in Section 4.

**2. Main result.** Throughout let $C(\mathbb{R}_0^+)$ be the space of continuous functions on $\mathbb{R}_0^+ = [0, \infty)$ endowed with the convergence of locally uniform convergence. The proof of the following theorem is found in Section 4.1.

THEOREM 2.1. *Suppose that* (1.1) *is satisfied. Let $L$ and $L_{|\cdot|}$ be slowly varying at infinity with $\lim_{t \to \infty} L(t) = \infty$ and*

$$\limsup_{t \to \infty} \frac{L_{|\cdot|}(t)}{L(t)} < \infty, \tag{2.1}$$

*where*

$$(2.2) \qquad L(t) = \int_0^t r^q(u)\,du, \qquad L_{|\cdot|}(t) = \int_0^t |r(u)|^q\,du, \qquad t > 0.$$

*Let $f \in L^2(\nu)$ with $Ef(N_0) = 0$ and Hermite rank $q = q(f) \geq 1$. Then the following assertions hold:*

(i) *In the sense of the finite dimensional distributions,*

$$(2.3) \qquad (2tL(t))^{-1/2} \int_0^{t\cdot} f(N_u)\,du \xrightarrow{d} \frac{c_q}{(q!)^{1/2}} W, \qquad t \to \infty.$$

(ii) *In addition, if there exists $R > 1$ such that $\sum_{k=q}^\infty \frac{|c_k|}{\sqrt{k!}} R^k$ is finite then* (2.3) *holds in $C(\mathbb{R}_0^+)$.*

REMARK 1. (i) It follows from Karamata's theorem, Fubini's theorem and the diagram formula (A.2) that

$$E\left[\int_0^t \mathcal{H}_q(N_s)\,ds\right]^2 = 2q! \int_0^t \int_0^v r^q(u)\,du\,dv \sim 2q!tL(t), \qquad t \to \infty,$$

and thus $L$ is eventually positive at infinity.

(ii) To set up weak convergence to a Brownian motion as in (2.3), it is well known that a possible norming function has to be regularly varying with index $1/2$ (cf. Lamperti [17] and the discussion in Taqqu [24]). Consequently, the assumption that $L$ is itself slowly varying is rather weak. Hariz [15] gives a functional version of his result under similar conditions on the Hermite expansion of $f$ and $N$. The existence of a representation of form (1.1) is equivalent to the absolute continuity of the spectral measure of $N$ with respect to Lebesgue measure (cf. Ibragimov and Linnik [16], Theorem 16.7.2). This covers a lot of important cases and, in particular, Ornstein–Uhenbeck processes driven by the fractional Brownian motion. Our proof of Theorem 1 relies on elegant results obtained by Peccati and Tudor [20] and formulated in terms of functionals on a Gaussian space. Therefore, we need this space controlled by a diffuse measure, which is taken to be the Lebesgue measure in the proof. It would be interesting to investigate whether the arguments by Peccati and Tudor [20] extend to other measures than the Lebesgue measure. However, this is beyond the scope of this paper.

(iii) For normal processes with a bounded spectral density, Chambers and Slud [11] have shown central limit theorems for more general functionals, including those depending on infinitely many coordinates. For the long-range dependence situation considered in our paper, the spectral density is usually unbounded, however. The result presented in this paper only covers functionals depending on one coordinate. A general theory has yet to be developed for general functionals of long-range dependent processes.

To prepare our analysis in the next section, we conclude this section with an application to fractional Ornstein–Uhlenbeck processes. Let $H \in (0,1)$. A fractional Brownian motion (FBM) $B^H = (B_t^H)_{t \in \mathbb{R}}$ with Hurst index $H$ is a centered Gaussian process with almost surely locally Hölder continuous paths of any order strictly smaller than $H$ and covariance function

$$(2.4) \qquad E B_t^H B_s^H = \tfrac{1}{2}(|t|^{2H} + |s|^{2H} - |t-s|^{2H}), \qquad s, t \in \mathbb{R}.$$

A choice $H = 1/2$ relates to a standard Brownian motion $W$ (cf. Bhattacharya and Waymire [3] for interpretations of Hurst index and Samorodnitsky and Taqqu [23] for further properties of FBM).

In particular, FBM is self-similar, that is, for all $c \in \mathbb{R}$, in the sense of finite dimensional distributions,

$$(2.5) \qquad (B_{ct}^H) \stackrel{d}{=} |c|^H (B_t^H).$$

Also, let $B_\gamma^H = (B_{\gamma,t}^H)_{t \geq 0}$ be the fractional Ornstein–Uhlenbeck process (FOUP), defined by

$$(2.6) \qquad B_{\gamma,t}^H = \int_0^t e^{-\gamma(t-s)} dB_s^H, \qquad t \geq 0, \ \gamma \in \mathbb{R},$$

where the integral in (2.6) converges pathwisely in the usual Riemann–Stieltjes sense (cf. Buchmann and Klüppelberg [9] and Cheridito, Kawaguchi and Maejima [12] for further properties of FOUP).

For $\gamma > 0$, the process $B_\gamma^H$ is asymptotically stationary. In the next corollary, we apply our result to the case when the Hermite rank $q = 2$. Up to some technicalities (cf. Section 4.2), it follows from the CLT by Hariz [15] and the NCLT by Taqqu [25] for $H \in (0, 3/4)$ and $H \in (3/4, 1)$, respectively. Theorem 2.1 deals with the remaining case $H = 3/4$. We conclude this section with a remark where we summarize the corresponding results for the general case and allow for arbitrary Hermite ranks $q$. We refer to Taqqu [25] and references therein for the definition and properties of Rosenblatt processes with Hurst index $H$. Throughout $\Gamma$ denotes Euler's Gamma function.

COROLLARY 2.1. *For all $H \in (0,1)$ and $\gamma > 0$, it holds in $C(\mathbb{R}_0^+)$*

$$(2.7) \quad g_H(t) \int_0^{t\cdot} \left[ (B_{\gamma,s}^H)^2 - \frac{\Gamma(2H+1)}{2\gamma^{2H}} \right] ds \stackrel{d}{\to} h_H(\gamma) \sigma_H X^H, \qquad t \to \infty,$$

*where, for $H \in (0, 3/4]$, $X^H$ is a Brownian motion and for $H \in (3/4, 1)$, $X^H$ is a Rosenblatt process with $E(X_1^H)^2 = 1$ and Hurst index $H$.*

*The quantities $g_H$, $\sigma_H$ and $h_H(\gamma)$ are given by the following formulas:*

$$(2.8) \quad g_H(t) = \begin{cases} t^{-1/2}, & H \in (0, 3/4), \\ (t \log t)^{-1/2}, & H = 3/4, \\ t^{1-2H}, & H \in (3/4, 1), \end{cases} \quad t > 1,$$

$$(2.9) \quad h_H(\gamma) = \begin{cases} \gamma^{-1/2-2H}, & H \in (0, 3/4), \\ \gamma^{-2}, & H \in [3/4, 1), \end{cases}$$

$$(2.10) \quad \sigma_H = \begin{cases} (2/\pi)^{1/2} \Gamma(2H+1) \\ \quad \times \sin(\pi H) \left( \int_0^\infty \dfrac{\xi^{2-4H}}{(1+\xi^2)^2} \, d\xi \right)^{1/2}, & H \in (0, 3/4), \\ 3/4, & H = 3/4, \\ H[(4H-2)/(4H-3)]^{1/2}, & H \in (3/4, 1). \end{cases}$$

REMARK 2. Let $q \in \mathbb{N}$, $H \in (0, 1)$ and $\gamma > 0$. In Corollary 2.1, we only deal with Hermite rank $q = 2$, which is sufficient for Section 3. It is natural to ask whether Corollary 2.1 can be extended to other Hermite polynomials. The answer is affirmative and as the arguments are very similar to the proof of Corollary 2.1, we only state the main results here without proof. To this end, let

$$(2.11) \quad \mu_H = \frac{2}{\Gamma(2H+1)}, \quad H \in (0, 1)$$

and observe that $E(B_{\gamma,t}^H)^2 \to \gamma^{-2H} \mu_H^{-1}$ (cf. Buchmann and Klüppelberg [9]).

In view of Lemma 4.1(iii) below, it follows from Taqqu [25] that, for $q \geq 1$ and $H > 1 - \frac{1}{2q}$, as $t \to \infty$, in $C(\mathbb{R}_0^+)$,

$$t^{q(1-H)-1} \int_0^{t \cdot} \mathcal{H}_q(\gamma^H \mu_H^{1/2} B_{\gamma,s}^H) \, ds$$

$$\xrightarrow{d} (2q!)^{1/2} [(2H-2)q+1][(2H-2)q+2]^{-1/2}$$

$$\times \left[ \frac{2H-1}{\Gamma(2H)} \right]^{q/2} \gamma^{q(H-1)} B^{q,H},$$

where $B^{q,H}$ is the Hermite process of order $q$ with Hurst index $H$ and unit variance $E(B_1^{q,H})^2 = 1$ (cf. [25] for the definition of higher order Hermite processes). Note that $B^{1,H} = B^H$ for $q = 1$ and $H \in (1/2, 1)$ whereas, for $q = 2$ and $H \in (3/4, 1)$, $B^{2,H} = X^H$ is the Rosenblatt process in Corollary 2.1.

If $q \geq 2$ and $H = 1 - 1/(2q)$, then it follows from Theorem 2.1 and Lemma 4.1(iii) that, as $t \to \infty$, in $C(\mathbb{R}_0^+)$,

$$(t \log t)^{-1/2} \int_0^{t \cdot} \mathcal{H}_q(\gamma^H \mu_H^{1/2} B_{\gamma,s}^H) \, ds \xrightarrow{d} (2q!)^{1/2} \left[ \frac{2H-1}{\Gamma(2H)} \right]^{q/2} \gamma^{q(H-1)} W.$$

If either, both $q = 1$ and $H \leq 1/2$, or, $q \geq 2$ and $H < 1 - 1/(2q)$, then we get from Hariz [15] and Lemma 4.1(iii) that, as $t \to \infty$, in $C(\mathbb{R}_0^+)$,

$$(2.12) \quad t^{-1/2} \int_0^{t\cdot} \mathcal{H}_q(\gamma^H \mu_H^{1/2} B_{\gamma,s}^H) \, ds \xrightarrow{d} \frac{(2q!)^{1/2}}{[2\Gamma(2H+1)]^{q/2}} \left(\frac{I_{q,H}}{\gamma}\right)^{1/2} W,$$

where

$$I_{q,H} = \int_0^\infty \left[\Gamma(2H+1)e^{-t} \right. $$
$$\left. + 2H\left[e^t \int_t^\infty e^{-s} s^{2H-1} \, ds - e^{-t} \int_0^t e^s s^{2H-1} \, ds\right]\right]^q dt.$$

Here, the integral converges absolutely, provided either $q = 1$ and $0 < H \leq 1/2$, or $q \geq 2$ and $0 < H < 1 - 1/(2q)$. Observe that $I_{q,1/2} = 2^q/q$ for all $q \in \mathbb{N}$. A more refined analysis shows that $I_{q,H} > 0$ for all $q \geq 2$ and $0 < H < 1 - 1/(2q)$.

However, note that $I_{1,H} = 0$ for all $0 < H < 1/2$ such that the limit in (2.12) is trivial in this case. To find a nontrivial limit, consider a function $\psi : [0, \infty) \to \mathbb{R}$ with $\psi(0) = 0$, that is, locally Hölder of order $\beta > 0$. For all $T > 0$ there, thus exists $C_T \in \mathbb{R}$ such that $|\psi(u) - \psi(v)| \leq C_T |u - v|^\beta$ for all $0 \leq u, v \leq T$. In this case,

$$\sup_{0 \leq v \leq T} \left| t \int_0^v e^{\gamma t(s-v)} \psi(s) \, ds - \frac{1}{\gamma} \psi(v) \right| \leq t^{-\beta} \frac{C_T}{\gamma^{1+\beta}} \left[\Gamma(\beta+1) + \sup_{0 \leq s < \infty} s^\beta e^{-s}\right],$$

for all $\gamma, t, T > 0$, and thus for all $\gamma > 0$, locally uniformly in $v \geq 0$,

$$(2.13) \quad t \int_0^v e^{\gamma t(s-v)} \psi(s) \, ds \to \frac{1}{\gamma} \psi(v), \quad t \to \infty.$$

Let $H \in (0,1)$, in view of (2.5), we get from simple substitutions and partial integrations that, for all $t > 0$,

$$(2.14) \quad t^{-H} \int_0^{t\cdot} B_{\gamma,s}^H \, ds = t^{1-H} \int_0^\cdot e^{\gamma t(s-\cdot)} B_{ts}^H \, ds \stackrel{d}{=} t \int_0^\cdot e^{\gamma t(s-\cdot)} B_s^H \, ds,$$

where the last identity holds in the sense of finite dimensional distributions. Clearly, we have $B_0^H = 0$ with probability one. Further, sample paths of $B^H$ are locally Hölder of any order strictly smaller than $H$, almost surely. Hence, the right-hand side in (2.14) converges to $B^H/\gamma$ locally uniformly, almost surely, as $t \to \infty$, by means of (2.13), and thus in $C(\mathbb{R}_0^+)$, as $t \to \infty$,

$$t^{-H} \int_0^{t\cdot} B_{\gamma,s}^H \, ds \xrightarrow{d} \gamma^{-1} B^H.$$

**3. The unit root problem.** We apply Corollary 2.1 to some functionals of FBM that occur in the unit-root problem in times series analysis. To illustrate this problem, let $\varepsilon = (\varepsilon_n)_{n\in\mathbb{N}}$ be a sequence of independent and identically distributed random variables with variance one and consider the first order autoregressive model,

$$X_n = \beta X_{n-1} + \varepsilon_n, \qquad n \in \mathbb{N}, \qquad X_0 = 0.$$

The parameter $\beta$ is unknown and has to be estimated from the observations $X_1, \ldots, X_n$. The least squares estimator for $\beta$ is given by the formula

$$\hat{b}_n = \hat{b}_n(X_0, \ldots, X_n) = \frac{\sum_{t=0}^{n-1} X_{t+1} X_t}{\sum_{t=0}^{n-1} X_t^2}.$$

For $|\beta| < 1$ (stationary regime), Mann and Wald [18] showed that

$$(3.1) \qquad \hat{\tau}_n = \left(\sum_{t=0}^{n-1} X_t^2\right)^{1/2} (\hat{b}_n - \beta) \xrightarrow{d} W_1, \qquad n \to \infty.$$

For the explosive case $|\beta| > 1$, (3.1) holds whenever $\varepsilon$ is a sequence of independent standard normal random variables (Anderson [1]).

On the other hand, (3.1) fails to hold for $\beta = 1$ even when $\varepsilon$ forms a sequence of independent standard normal random variables. White [26] and Rao [22] showed that it is a functional of Brownian motion $W = (W_t)_{0 \le t \le 1}$, that is,

$$(3.2) \qquad \hat{\tau}_n \xrightarrow{d} \bar{\tau}(0) = \tfrac{1}{2}[W_1^2 - 1]\left[\int_0^1 W_s^2 \, ds\right]^{-1/2}, \qquad n \to \infty.$$

This contrast between (3.1) and (3.2) is a typical example of a critical phenomenon. The parameter value $\beta = 1$ comprises a singularity and there is a lack of smooth transition of the limiting distribution of $\hat{\tau}_n$ when $\beta$ is close to one. For finite sample analysis or tests under local alternatives, a key question becomes that if $\beta$ is close to one, what kind of approximation should be used for $\hat{\tau}_n$? An answer to this question is given in Chan and Wei [10], where the following class of models is proposed:

$$(3.3) \qquad X_t^{(n)} = \beta_n X_{t-1}^{(n)} + \varepsilon_t, \qquad n \in \mathbb{N}, \qquad X_0^{(n)} = 0.$$

Suppose that there exists $\gamma \in \mathbb{R}$ such that $\beta_n = 1 - \gamma/n$. Then it is shown by Chan and Wei [10] that

$$(3.4) \qquad \hat{\tau}_n \xrightarrow{d} \bar{\tau}(\gamma) = \frac{\int_0^1 W_{\gamma,s} \, dW_s}{\sqrt{\int_0^1 W_{\gamma,s}^2 \, ds}}, \qquad n \to \infty,$$

where $W_\gamma = B_\gamma^{1/2}$ is the Ornstein–Uhlenbeck process driven by Brownian motion and the integral on the right-hand side is defined in the Itô's sense.

In [8], a generalized asymptotic theory of ordinary least squares estimators for a large class of possibly strongly dependent noise sequences was given. The result covers certain stationary and ergodic sequences $\varepsilon = (\varepsilon_n)_{n \in \mathbb{N}}$ with mean zero and finite variance $\sigma^2 = E\varepsilon_2^2$, where the associated partial sum process is in the domain of attraction of $B^H$ for some $H \in (0,1)$. In this case, there exist matrices $D_n = D_n(H) \in \mathbb{R}^{2 \times 2}$, $n \in \mathbb{N}$, $H \in (0,1)$, such that, as $n \to \infty$,

$$(3.5) \quad D_n \begin{pmatrix} \hat{\tau}_n \\ \hat{b}_n - \beta_n \end{pmatrix} \overset{d}{\to} 1_{H \geq 1/2} \begin{pmatrix} \tau_{H,3}(\gamma) \\ \tau_{H,4}(\gamma) \end{pmatrix} - 1_{H \leq 1/2} \frac{\sigma^2}{2} \begin{pmatrix} \tau_{H,1}(\gamma) \\ \tau_{H,2}(\gamma) \end{pmatrix},$$

where $1_A$ denotes the indicator function of some set $A$. Further, $\tau_{H,i}(\gamma)$ denotes the $i$th component of the vector $\tau_H(\gamma) \in \mathbb{R}^4$ defined by

$$(3.6) \qquad \tau_H(\gamma) = \left( \Theta_H(\gamma)', \left( \frac{(B^H_{\gamma,1})^2}{2} + \gamma \int_0^1 (B^H_{\gamma,s})^2\, ds \right) \Theta_H(\gamma)' \right)',$$

$$(3.7) \qquad \Theta_H(\gamma) = \begin{pmatrix} \left( \int_0^1 (B^H_{\gamma,s})^2\, ds \right)^{-1/2} \\ \left( \int_0^1 (B^H_{\gamma,s})^2\, ds \right)^{-1} \end{pmatrix},$$

where $x'$ denotes the transpose of a vector $x \in \mathbb{R}^n$.

It follows from Itô's formula that

$$(3.8) \qquad \bar{\tau}(\gamma) = \tau_{1/2,3}(\gamma) - \tfrac{1}{2}\tau_{1/2,1}(\gamma),$$

and thus the limit in (3.4) can be derived from (3.5) and (3.8), provided $\sigma^2 = 1$ (cf. [8] for details).

It was shown in [10], that $\bar{\tau}(\gamma)$ forms a continuous family indexed by $\gamma \in [-\infty, \infty]$, in particular, we have with (3.1) in view,

$$(3.9) \qquad \bar{\tau}(\gamma) \overset{d}{\to} W_1, \qquad |\gamma| \to \infty.$$

Next we want to extend (3.9) to all $H \in (0,1)$. Therefore, we consider the vector $\tau_H(\gamma)$ and its asymptotic properties for $|\gamma| \to \infty$. In Remark 5 below, we show how the limit in (3.9) is derived from the general theory.

There are two cases, that is, $\gamma \to \infty$ or $\gamma \to -\infty$. Both cases exhibit quite different qualitative behavior. With (3.9) in mind, this seems to be puzzling at the first glance. However, the underlying dynamics are rather different as formally taking $\gamma$ to $+\infty$ or $-\infty$ in (3.3) corresponds to an infinitesimal return to the stationary regime or the explosive one.

We first tackle $\gamma \to \infty$. The corresponding theorem is implied by Corollary 2.1 as it follows from the self-similarity of FBM that we may rewrite $\tau_H(\gamma)$ into a functional of integrated squares of $B^H_1$. In view of Corollary 2.1 for Hermite rank $q = 2$, we thus expect two types of limit distributions depending on $H \leq 3/4$ or $H > 3/4$ with an additional logarithm in the norming functions for $H = 3/4$ (cf. Section 4.3 for a proof).

THEOREM 3.1. *Let $(Z,Y)$ a standard normal random vector. For $H \in (3/4,1)$, let $R_H = R^H(1)$ be a Rosenblatt distributed random variable with $ER_H^2 = 1$, independent of $Y$.*

(i) *For $H \in (0, 3/4]$ and $\gamma > 1$, there exist $\Sigma_H \in \mathbb{R}^{4\times 2}$, $b_H(\gamma) \in \mathbb{R}^4$ and $D_H(\gamma) \in \mathbb{R}^{4\times 4}$ such that*

$$D_H(\gamma)(\tau_H(\gamma) - b_H(\gamma)) \xrightarrow{d} \Sigma_H \begin{pmatrix} Z \\ Y^2 \end{pmatrix}, \qquad \gamma \to \infty.$$

(ii) *For $H \in (3/4, 1)$ and $\gamma > 1$ there exist $\Sigma_H \in \mathbb{R}^{4\times 2}$, $b_H(\gamma) \in \mathbb{R}^4$ and $D_H(\gamma) \in \mathbb{R}^{4\times 4}$ such that*

$$D_H(\gamma)(\tau_H(\gamma) - b_H(\gamma)) \xrightarrow{d} \Sigma_H \begin{pmatrix} R_H \\ Y^2 \end{pmatrix}, \qquad \gamma \to \infty.$$

REMARK 3. For $\gamma > 1$ and $H \in (0,1)$, the scaling matrices $D_H(\gamma)$ have the following forms:

$$D_H(\gamma) = \begin{cases} \operatorname{diag}(\gamma^{1/2-H}, \gamma^{1/2-2H}, \gamma^{H-1/2}, 1), & H \in (0, 3/4), \\ \operatorname{diag}(\gamma^{-1/4}(\log\gamma)^{-1/2}, \gamma^{-1}(\log\gamma)^{-1/2}, \\ \qquad \gamma^{1/4}(\log\gamma)^{-1/2}, 1), & H = 3/4, \\ \operatorname{diag}(\gamma^{2-3H}, \gamma^{2-4H}, \gamma^{1-H}, 1), & H \in (3/4, 1). \end{cases}$$

For $\gamma > 0$ and $H \in (0,1)$, the matrices $\Sigma_H$ and centering vectors $\xi_H(\gamma)$ are given by the following formulas

$$\Sigma_H = \begin{pmatrix} -\kappa_H \mu_H & 0 \\ -2\kappa_H \mu_H^{3/2} & 0 \\ \kappa_H & 0 \\ 0 & 1/2 \end{pmatrix}, \qquad b_H(\gamma) = \begin{pmatrix} \mu_H^{1/2} \gamma^H \\ \mu_H \gamma^{2H} \\ \mu_H^{-1/2} \gamma^{1-H} \\ \gamma \end{pmatrix},$$

where $\mu_H$ is defined in (2.11) and we set

$$\kappa_H = \begin{cases} \pi^{-1/2} \sin(\pi H) \Gamma(2H+1)^{1/2} \left( \int_0^\infty \dfrac{\xi^{2-4H}}{(1+\xi^2)^2}\, d\xi \right)^{1/2}, & H \in (0, 3/4), \\ (3/8)^{1/2} \pi^{-1/4}, & H = 3/4, \\ 2^{-1/2} \left( \dfrac{H}{4H-3} \right)^{1/2} \Gamma(2H-1)^{-1/2}, & H \in (3/4, 1). \end{cases}$$

We conclude this section with a corresponding result dealing with the other case, namely $\gamma \to -\infty$, together with some remarks (cf. Section 4.4 for its proof).

THEOREM 3.2. *Let $H \in (0, 1)$. Then for $\gamma \to -\infty$,*

$$\operatorname{diag}(|\gamma|^{-(2H+1)/2} e^{|\gamma|}, |\gamma|^{-2H-1} e^{2|\gamma|}, |\gamma|^{(2H-1)/2}, |\gamma|^{-1} e^{|\gamma|}) \tau_H(\gamma)$$

$$\overset{d}{\to} \mathrm{diag}\left(\frac{2}{\Gamma(2H+1)^{1/2}}, \frac{4}{\Gamma(2H+1)}, \Gamma(2H+1)^{1/2}, 2\right)$$
$$\times (|Z|^{-1}, Z^{-2}, Y\mathrm{sign}(Z), Y/Z)',$$

where $(Y, Z)$ is standard normal and

$$\mathrm{sign}(x) = 1, \quad x \geq 0, \quad \mathrm{sign}(x) = -1, \quad x < 0.$$

REMARK 4. In Theorem 3.2, observe that $Y/Z$ is a Cauchy random variable. A similar result was obtained by Anderson [1] for the ordinary least squares estimator $\hat{b}_n$ for nonstationary AR(1) models, that is, $|\beta| > 1$, with independent standard Gaussian innovations.

REMARK 5. As indicated, we aim for recovering the limit in (3.9) by Theorems 3.1–3.2. An easy computation yields $\mu_{1/2} = 2$ and $\kappa_{1/2} = 1/2$, by virtue of the identity $\int_0^\infty d\xi/(1+\xi^2)^2 = \pi/4$. By Theorem 3.1,

$$\left(\begin{pmatrix}\tau_{1/2,1}(\gamma)\\ \tau_{1/2,3}(\gamma)\end{pmatrix} - \gamma^{1/2}\begin{pmatrix}2^{1/2}\\ 2^{-1/2}\end{pmatrix}\right) \overset{d}{\to} Z\begin{pmatrix}-1\\ 1/2\end{pmatrix}, \quad \gamma \to \infty$$

and thus by (3.8),

$$\bar{\tau}(\gamma) = (\tau_{1/2,3}(\gamma) - 2^{-1/2}\gamma^{1/2}) - \tfrac{1}{2}(\tau_{1/2,1} - 2^{1/2}\gamma^{1/2}) \overset{d}{\to} Z, \quad \gamma \to \infty.$$

On the other hand, Theorem 3.2 yields

$$\begin{pmatrix}|\gamma|^{-1}e^{|\gamma|}\tau_{1/2,1}(\gamma)\\ \tau_{1/2,3}(\gamma)\end{pmatrix} \overset{d}{\to} \begin{pmatrix}2|Z|^{-1}\\ Y\mathrm{sign}(Z)\end{pmatrix}, \quad \gamma \to -\infty.$$

By means of (3.8), we get, for $\gamma \to -\infty$,

$$\bar{\tau}(\gamma) = \tau_{1/2,3}(\gamma) - \tfrac{1}{2}|\gamma|e^{-|\gamma|}(|\gamma|^{-1}e^{|\gamma|}\tau_{1/2,1}(\gamma)) \overset{d}{\to} Y\mathrm{sign}(Z) \overset{d}{=} Z.$$

**4. Proofs.**

4.1. *Proof of Theorem 2.1.* (i) First, we wish to establish (i) for the special choice $f = \mathcal{H}_q$, $q \in \mathbb{N}$. We make use of a CLT for vector-valued multiple stochastic integrals by Peccatti and Tudor [20] (cf. [19], Section 1.1, for the necessary background). Secondly, we extend the result by a reduction, this being similar as in [24].

To this end, fix $q \in \mathbb{N}$ and let

$$A_{q,t}(v_1, v_2) = (2tL(t))^{-1/2}\int_{v_1 t}^{v_2 t}\mathcal{H}_q(N_u)\,du, \quad v_2 \geq v_1 \geq 0.$$

In order to show (2.3) for the finite dimensional distributions for $f = \mathcal{H}_q$, it suffices to show that, as $t \to \infty$,

(4.1) $\quad (A_{q,t}(v_{2j-1}, v_{2j}))_{1 \leq j \leq d} \overset{d}{\to} (q!)^{1/2}(W_{v_{2j}} - W_{v_{2j-1}})_{1 \leq j \leq d},$

for all $d \in \mathbb{N}$ and $v_{2d} > v_{2d-1} \geq v_{2d} > \cdots \geq v_2 > v_1 \geq 0$.

Fix $d \in \mathbb{N}$ and $v_{2d} > v_{2d-1} \geq v_{2d} > \cdots \geq v_2 > v_1 \geq 0$. By Proposition 1.1.4 in [19] and a stochastic version of the stochastic Fubini's theorem, it follows from (2.2) that, almost surely for $1 \leq j \leq d$,

$$
(4.2) \quad \begin{aligned} A_{q,t}(v_{2j-1}, v_{2j}) \\ = (2tL(t))^{-1/2} \int_{\mathbb{R}^q} \int_{v_{2j-1}t}^{v_{2j}t} \prod_{i=1}^{q} a(t-u_i) \, dW^{\otimes q}(u_1, \ldots, u_d), \end{aligned}
$$

where the right-hand side is understood in terms of multiple Wiener integrals (note that the Hermite polynomials in (1.2) have a different leading coefficient than the ones in [19]).

Suppose that we can show that, for all $v_4 > v_3 \geq v_2 > v_1 \geq 0$, as $t \to \infty$,

$$(4.3) \quad EA_{q,t}^2(v_1, v_2) \to q! E(W_{v_2} - W_{v_1})^2,$$

$$(4.4) \quad E[A_{q,t}(v_1, v_2) A_{q,t}(v_3, v_4)] \to 0$$

and

$$(4.5) \quad EA_{q,t}^4(v_1, v_2) \to q! E(W_{v_2} - W_{v_1})^4.$$

Then (4.1) is immediate in view of (4.2) and Theorem 1 in [20]. To show (4.3), we make use of the diagram formula, which together with some necessary graph-theoretical concepts, are stated in the Appendix.

Fix $v_2 > v_1 \geq 0$. As $r$ is symmetric and $L$ is slowly varying, we get from Fubini's theorem and Karamata's theorem that, as $t \to \infty$,

$$(4.6) \quad \int_{tv_1}^{tv_2} \int_{tv_1}^{tv_2} r^q(u_1 - u_2) \, du_1 \, du_2 = 2 \int_0^{t(v_2-v_1)} L(u) \, du \sim 2(v_2 - v_1) tL(t),$$

and thus by (A.2), as $t \to \infty$,

$$EA_{q,t}^2(v_1, v_2) = 2q!(2tL(t))^{-1} \int_{tv_1}^{tv_2} \int_{tu}^{tw} r^q(u_1 - u_2) \, du_1 \, du_2 \to q!(v_2 - v_1).$$

This completes the proof of (4.3).

Next pick $v_4 > v_3 \geq v_2 > v_1 \geq 0$. By simple substitutions, observe that

$$\begin{aligned} \int_{tv_3}^{tv_4} \int_{tv_1}^{tv_2} r^q(u_1 - u_2) \, du_1 \, du_2 \\ = \left[ \int_0^{t(v_4-v_1)} + \int_0^{t(v_3-v_2)} - \int_0^{t(v_4-v_2)} - \int_0^{t(v_3-v_1)} \right] L(u) \, du, \end{aligned}$$

for $t > 0$, and thus by means of Karamata's theorem, as $t \to \infty$,

$$(4.7) \quad \int_{tv_3}^{tv_4} \int_{tv_1}^{tv_2} r^q(u_1 - u_2) \, du_1 \, du_2 = o(tL(t)), \qquad t \to \infty.$$

Now (4.4) is implied by (A.2) and (4.7) since, as $t \to \infty$,

$$E[A_{q,t}(v_1,v_2)A_{q,t}(v_3,v_4)]$$
$$= (2q!)(2tL(t))^{-1} \int_{tv_3}^{tv_4} \int_{tv_1}^{tv_2} r^q(u_1-u_2)\, du_1\, du_2 \to 0.$$

To show (4.5) let $\varepsilon > 0$ and recall (2.2). It follows from (2.1) and $\lim_{t \to \infty} L(t) = \infty$ that $I := \int_0^\infty |r(u)|^q\, du$ is infinite. By assumption, $L_{|\cdot|}$ is also slowly varying. Consequently, by Hölder's inequality and Karamata's theorem, for all $v > \varepsilon > 0$,

$$\int_0^{tv} |r(u)|^m\, du \leq (\varepsilon t)^{1-m/q} L_{|\cdot|}^{m/q}(\varepsilon t) + [(v-\varepsilon)t]^{1-m/q}[L_{|\cdot|}(vt) - L_{|\cdot|}(\varepsilon t)]^{m/q}$$
$$\sim (\varepsilon t)^{1-m/q} L_{|\cdot|}^{m/q}(t), \qquad t \to \infty,\ 1 \leq m < q,$$

and thus by (2.1),

$$(4.8) \qquad \int_0^{tv} |r(u)|^m\, du = o(t^{1-m/q} L(t)^{m/q}), \qquad t \to \infty,\ 1 \leq m < q,\ v > 0.$$

On the other hand, it follows from the inequality between geometric and arithmetic mean and simple substitutions that for all $m \in \mathbb{N}$, $v > 0$ and $0 \leq w_1, \ldots, w_m \leq vt$,

$$(4.9) \qquad \int_0^{vt} \prod_{i=1}^m |r(u-w_i)|\, du \leq 2 \int_0^{vt} |r(u)|^m\, du, \qquad t > 0.$$

We return to the proof of (4.5). Fix $v_2 > v_1 \geq 0$. Using the notation of the Appendix, let $\mathcal{D}(4,q)$ be the set of all diagrams with 4 levels and rows of length $q$. Let $\mathcal{D}_q^r$ be the subset of all $D \in \mathcal{D}(4,q)$ such that the following implication holds for all $e, f \in \mathcal{V}(D)$:

$$m(e) = m(f) \quad \Rightarrow \quad M(e) = M(f).$$

If $D \in \mathcal{D}(4,q) \setminus \mathcal{D}_q^r$ then, for $1 \leq i \leq 3$ and $i < j \leq 4$, there exist integers $0 \leq m_{i,j} \leq q$, additionally satisfying $m_{1,1} + m_{1,2} + m_{1,3} + m_{1,4} = q$ and $1 \leq m_{2,3} + m_{2,4} < q$, such that

$$\left| \int_{tv_1}^{tv_2} \int_{tv_1}^{tv_2} \int_{tv_1}^{tv_2} \int_{tv_1}^{tv_2} \prod_{e \in \mathcal{E}(D)} r(u_{m(e)} - u_{M(e)})\, du_1\, du_2\, du_3\, du_4 \right|$$
$$\leq \int_0^{tv_2} \int_0^{tv_2} \int_0^{tv_2} \int_0^{tv_2} \prod_{i=1}^3 \prod_{j=i+1}^4 |r(u_i - u_j)|^{m_{i,j}}\, du_1\, du_2\, du_3\, du_4.$$

Thus, by (4.8)–(4.9) and (2.1), as $t \to \infty$, we have

$$\left| \int_{tv_1}^{tv_2} \int_{tv_1}^{tv_2} \int_{tv_1}^{tv_2} \int_{tv_1}^{tv_2} \prod_{e \in \mathcal{E}(D)} r(u_{m(e)} - u_{M(e)}) \, du_1 \, du_2 \, du_3 \, du_4 \right|$$

$$\leq 8u_2 t L_{|\cdot|}(v_2 t) \int_0^{tv_2} |r(u)|^{m_{3,4}} \, du_1 \int_0^{tv_2} |r(u)|^{m_{2,3}+m_{2,4}} \, du_1$$

$$= o(t^2 L^2(t)), \qquad t \to \infty.$$

Consequently, we get from (A.2) and Fubini's theorem that, as $t \to \infty$,

$$E A_{q,t}^4(v_1, v_2) = 3(q!)^2 \left( \int_{tv_1}^{tv_2} \int_{tv_1}^{tv_2} r^q(u_1 - u_2) \, du_1 \, du_2 \right)^2 + o(t^2 L^2(t)).$$

This completes the proof of (4.5) in view of (4.6).

To summarize, the assertions (4.3)–(4.5) are all in place, and thus (4.1) is implied by (4.2) and Theorem 1 in [20]. Further, this completes the proof of (i) for the choice $f = \mathcal{H}_q$. To extend (2.3) to all $f \in L^2(\nu)$, recall that $\lim_{t \to \infty} L(t) = \infty$, and thus $\lim_{t \to \infty} L_{|\cdot|}(t) = \infty$ by (2.1). Further, it follows from (2.2) that a spectral density exists (cf. Ibragimov and Linnik [16], Theorem 16.7.2). In particular, $r$ must be continuous and we have $\lim_{t \to \infty} r(t) = 0$ by means of the Riemann–Lebesgue lemma. If $\int_0^\infty |r(u)|^{q+1} \, du$ is infinite, then we get from l'Hôspital's rule that, for all $v > 0$,

$$(4.10) \qquad \lim_{t \to \infty} \int_0^{vt} |r(u)|^{q+1} \, du / L_{|\cdot|}(t) = \lim_{t \to \infty} |r(t)| = 0.$$

Otherwise, if $\int_0^\infty |r(u)|^{q+1} \, du$ is finite, (4.10) is trivial.

On the other hand, it follows from stationarity and Cauchy–Schwarz inequality that $|r(t)| \leq t$ for all $t \in \mathbb{R}$, and thus for all $v, t > 0$ and $m \geq q + 1$,

$$\int_0^{vt} \int_0^{vt} |r(u_1 - u_2)|^m \, du_1 \, du_2 \leq 2tv \int_0^{vt} |r(u)|^{q+1} \, du.$$

Hence, we get from (2.1) and (4.10) that, for all $v > 0$, as $t \to \infty$,

$$(4.11) \quad M_t(v) := \sup_{m > q} \int_0^{vt} \int_0^{vt} |r(u_1 - u_2)|^m \, du_1 \, du_2 = o(tL(t)), \qquad t \to \infty.$$

Now pick $f \in L^2(\nu)$ with $Ef(N_0) = 0$ and Hermite rank $q(f) = q$; by (A.2), for all $t > 0$, $d \in \mathbb{N}$, $\theta_1, \ldots, \theta_d \in \mathbb{R}$ and $v_d > v_{d-1} > \cdots > v_1 > v_0 = 0$,

$$E \left( \sum_{k=1}^d \theta_i \int_{v_{k-1}t}^{v_k t} f(N_u) - \frac{c_q}{q!} H_q(N_u) \, du \right)^2$$

$$\leq \sum_{1 \leq i_1, i_2 \leq d} |\theta_{i_1} \theta_{i_2}| \sum_{k=q+1}^\infty \frac{c_k^2}{k!} \int_0^{v_d t} \int_0^{v_d t} |r(u_1 - u_2)|^m \, du_1 \, du_2$$

$$\leq M_t(v_d) \sum_{1 \leq i_1, i_2 \leq d} |\theta_{i_1} \theta_{i_2}| \sum_{k=q+1}^{\infty} \frac{c_k^2}{k!},$$

where $c_k$ are the quantity defined in (1.4). Note that the sum on the right-hand side is finite since $f \in L^2(\nu)$. We have already shown that (i) holds for the particular choice $f = \mathcal{H}_q$. In view of (4.11) and the chain of inequalities in the last display, the proof of assertion (i) is completed by the Cramér–Wold device.

(ii) Let $f \in L^2(\nu)$ with $Ef(N_u) = 0$ and Hermite rank $q(f) = q$ such that there exists $R > 1$ with

(4.12) $$\sum_{k=q}^{\infty} \frac{|c_k|}{\sqrt{k!}} R^k < \infty.$$

In view of (i), it remains to show that $(2tL(t))^{-1/2} \int_0^{t \cdot} f(N_u)$ is tight in $C(\mathbb{R}_0^+)$, as $t \to \infty$. In view of Billingsley's tightness condition ([6], Theorem 13.14) and the Cauchy–Schwarz inequality, for $x_0 > 0$ it suffices to find $\theta > 0$, $t_0 > 0$ and $C$ such that, for all $0 \leq x \leq x_0$ and $t > t_0$,

(4.13) $$E\left( (2tL(t))^{-1/2} \int_0^{xt} f(N_u) \, du \right)^{2(1+\theta)} \leq C x^{1+\theta}.$$

Therefore, recall the following inequality as shown by Hariz [15] (cf. proof of Theorem 1(ii) in [15]): for all $0 < \theta \leq 1$, $t > 0$,

$$E\left( \int_0^t f(N_u) \, du \right)^{2(1+\theta)} \leq \left[ \sqrt{2t} \sum_{m=q}^{\infty} \frac{3^{m\theta/(1+\theta)} |c_m|}{\sqrt{m!}} \left( \int_0^t |r(u)|^m \, du \right)^{1/2} \right]^{2(1+\theta)}.$$

Consequently, for all $0 < \theta \leq 1$ and $t > 0$, we have by stationarity and Cauchy–Schwarz inequality,

$$E\left( \int_0^t f(N_u) \, du \right)^{2(1+\theta)} \leq \left( \sqrt{2tL_{|\cdot|}(t)} \sum_{m=q}^{\infty} \frac{3^{m\theta/(1+\theta)} |c_m|}{\sqrt{m!}} \right)^{2(1+\theta)}.$$

Thus, it follows from (4.12) that there exist $\theta \in (0, 1]$ and a constant $C'$ such that, for all $t > 0$,

$$E\left( \int_0^t f(N_u) \, du \right)^{2(1+\theta)} \leq C'(2tL_{|\cdot|}(t))^{1+\theta}.$$

Consequently, for all $x_0 > 0$ and $0 \leq x \leq x_0$, we have

$$E\left( (2tL(t))^{-1/2} \int_0^{xt} f(N_u) \, du \right)^{2(1+\theta)} \leq C'\left( \frac{L_{|\cdot|}(x_0 t)}{L(t)} \right)^{1+\theta} x^{1+\theta}.$$

As $L_{|\cdot|}$ is slowly varying, it follows from (2.1) that there is indeed a constant $C \in \mathbb{R}^+$ and $t_0 > 0$ such that (4.13) holds for all $t > 0$ and $0 \leq x \leq x_0$, giving (ii).

4.2. *Proof of Corollary 2.1.* For $H \in (0,1)$, we set $\mu_H = 2/\Gamma(2H+1)$ [cf. (2.11)]. With the help of $\mu_H$, we define a normal process $N_\gamma^H = (N_{\gamma,t}^H)_{t \in \mathbb{R}}$ by

$$(4.14) \quad N_{\gamma,t}^H = \gamma^H \mu_H^{1/2} \int_{-\infty}^t e^{-\gamma(t-s)} \, dB_s^H, \qquad t \in \mathbb{R}, \ H \in (0,1), \ \gamma > 0,$$

and set $r_{H,\gamma}(t) = E N_{\gamma,t}^H N_{\gamma,0}^H$, $t \in \mathbb{R}$. Then we have the following lemma:

LEMMA 4.1. *Let $\gamma > 0$, $H \in (0,1)$ and $q \in \mathbb{N}$.*
(i) *If $H \neq 1/2$, then*

$$r_{H,\gamma}(t) = \gamma^{2H-2} \frac{2H-1}{\Gamma(2H)} t^{2H-2} + O(t^{2H-4}), \qquad t \to \infty,$$

*whereas $\rho_{1/2,\gamma}(t) = e^{-\gamma |t|}$ for all $t \in \mathbb{R}$.*
(ii)

$$\int_\mathbb{R} r_{H,\gamma}^2(t) \, dt = \gamma^{-1} \sin^2(\pi H) \frac{4}{\pi} \int_0^\infty \frac{\xi^{2-4H}}{(1+\xi^2)^2} \, d\xi.$$

(iii)

$$(4.15) \qquad \int_0^\infty |\mathcal{H}_q(\gamma^H \mu_H^{1/2} B_{\gamma,s}^H) - \mathcal{H}_q(N_{\gamma,s}^H)| \, ds < \infty \qquad a.s.$$

PROOF. (i) For $H = 1/2$, this is a well-known property of the classical Ornstein–Uhlenbeck process driven by the Wiener process. For $H \neq 1/2$, assertion (i) is found in Theorem 2.3 in [12].

(ii) The spectral density $f_H$ of $N_\gamma^H$, that is, $r_{H,\gamma}(t) = \int_\mathbb{R} e^{it\xi} f_{H,\gamma}(\xi) \, d\xi$, has the following form ([12], (2.2)):

$$(4.16) \quad f_{H,\gamma}(\xi) = \gamma^{2H} \pi^{-1} \sin(\pi H) \frac{|\xi|^{1-2H}}{\gamma^2 + \xi^2}, \qquad \xi \in \mathbb{R}, \ H \in (0,1).$$

In particular, $\int_\mathbb{R} f_{H,\gamma}^2(\xi) \, d\xi$ is finite if and only if $H \in (0, 3/4)$. In this case, the identity in (ii) follows from Plancherel's theorem and simple substitutions.

(iii) Fix $q \in \mathbb{N}$ and $H \in (0,1)$. Then we find a constant $C_q \in (0, \infty)$ such that $|\mathcal{H}_q(x) - \mathcal{H}_q(y)| \leq C_q |x-y|(1 + |x|^q + |y|^q)$ for all $x, y \in \mathbb{R}$, where $\mathcal{H}_q$ is defined in (1.2). Note that almost surely,

$$(4.17) \qquad \gamma^H \mu_H B_{\gamma,t}^H = N_{\gamma,t}^H - e^{-\gamma t} N_{\gamma,0}^H, \qquad t \geq 0,$$

which completes the proof of (iii) as almost surely,

$$\int_0^\infty |\mathcal{H}_q(\gamma^H \mu_H B_{\gamma,s}^H)^2 - \mathcal{H}_q(N_{\gamma,s}^H)| \, ds$$

$$\leq C_q |N_{\gamma,0}^H| \int_0^\infty e^{-\gamma u} (1 + 2^{q+1} |N_{\gamma,s}^H|^q + 2^q |N_{\gamma,0}^H|^q) \, ds < \infty.$$

We return to the proof of Corollary 2.1: If $H \in (0, 3/4)$ then (2.7) follows from the functional CLT by Hariz [15] and Lemma 4.1. By employing the machinery developed by Pipiras and Taqqu ([21], Theorem 3.2), it is straightforward to show that $N_\gamma^H$ satisfies the assumptions imposed by Taqqu in [25] for all $H \in (1/2, 1)$, and thus (2.7) is in place for all $H \in (3/4, 1)$. Recall that the spectral measure of $N_\gamma^H$ is absolutely continuous with respect to Lebesgue measure [cf. (4.16)], and thus FOUP admits the representation in (2.2) (cf. Ibragimov and Linnik [16], Theorem 16.7.2). For $H = 3/4$, (2.7) is thus implied by Theorem 2.1 and Lemma 4.1(i). □

4.3. *Proof of Theorem 3.1.* In order to show Theorem 3.1, we need the following modification of Corollary 2.1.

LEMMA 4.2. *For $H \in (0, 3/4]$ let $X^H$ be standard normal random variable and, for $H \in (3/4, 1)$, let $X^H$ be a Rosenblatt distributed random variable with unit variance. Then, for all $H \in (0, 1)$, we have*

$$(4.18) \quad \begin{pmatrix} g_H(t) \int_0^t [(B_{1,s}^H)^2 - \mu_H^{-1}] \, ds \\ (B_{1,t}^H)^2 \end{pmatrix} \xrightarrow{d} \begin{pmatrix} \sigma_H X^H \\ \mu_H^{-1} Y^2 \end{pmatrix}, \qquad t \to \infty,$$

*where $Y$ is a standard normal random variable, independent of $X^H$, and the quantities $g_H$ and $\sigma_H$ and $\mu_H$ are defined in (2.8) and (2.10) and (2.11), respectively.*

PROOF. Let $H \in (0, 1)$ and recall (2.11). It is clear that $(B_{1,t}^H)^2 \xrightarrow{d} \mu_H^{-1} Y^2$, as $t \to \infty$. Hence, weak convergence of the marginals follows from Corollary 2.1 and Lemma 4.1(ii) and (iii). It remains to show that the marginals are asymptotically independent. This is implied by (4.17) and Lemma 4.1(iv), provided we can show that the components of the vector in the following display are asymptotically independent, as $t \to \infty$:

$$(4.19) \quad \begin{pmatrix} g_H(t) \int_0^t \mathcal{H}_2(N_{1,u}^H) \, du \\ N_{1,t}^H \end{pmatrix},$$

where $N_{1,\cdot}^H$ is defined in (4.14). In view of the stationarity of $N_1^H$ and Lemma 4.1(iii), it suffices to show that, for all $\theta_1, \theta_2 \in \mathbb{R}$,

$$(4.20) \quad \lim_{t \to \infty} E\left[ (e^{i\theta_1 N_{1,0}^H} - Ee^{i\theta_1 N_{1,0}^H}) \exp\left( i\theta_2 g_H(t) \int_{-t}^0 \mathcal{H}_2(N_{1,u}^H) \, du \right) \right] = 0.$$

For $n \in \mathbb{N}$, let $M_n = E[e^{i\theta_2 N_{1,0}^H} | \mathcal{F}_n]$ where $\mathcal{F}_n$ is the $\sigma$-algebra generated by $(N_{1,s}^H)_{s \leq -n}$. Then $(M_n, \mathcal{F}_n)_{n \in \mathbb{N}}$ is a (bounded) backward martingale. By virtue of the convergence theorem for backward martingales (cf. Hall and Heyde [14], Theorem 2.6), note that $\lim_{n\to\infty} M_n = E\exp(i\theta_2 N_{1,0}^H)$ almost surely.

For $\gamma > 1$ pick $n_t \in \mathbb{N}$, $n_t \leq t$ such that $n_t \to \infty$ and $g_H(t)/g_H(n_t) = 0$ as $t \to \infty$. In view of the triangular inequality, conditioning on $\mathcal{F}_{n_t}$, this yields the following inequality, for $t > 1$,

$$\left| E\left[(e^{i\theta_1 N_{1,0}^H} - Ee^{i\theta_1 N_{1,0}^H})\exp\left(i\theta_2 g_H(t)\int_{-t}^0 \mathcal{H}_2(N_{1,u}^H)\,du\right)\right]\right|$$
$$\leq 2E\left|\exp\left(i\theta_2 g_H(t)\int_{-n_t}^0 \mathcal{H}_2(N_{1,u}^H)\,du\right) - 1\right| + E|M_{n_t} - E\exp(i\theta_1 N_{1,0}^H)|.$$

By the choice of $n_t$, $g_H(\gamma)\int_{n_t}^0 \mathcal{H}_2(N_u^H)$ tends to zero in probability for $t \to \infty$. Hence, the right-hand side in the last display tends to zero as $t \to \infty$, giving (4.20). □

We return to the proof of Theorem 3.1. Fix $H \in (0,1)$, for $\gamma > 0$, it follows from the self-similarity of $B^H$ [cf. (2.5)] that, for the finite dimensional distributions we have

$$(4.21) \quad (B_{\gamma,t}^H)_{t\geq 0} = \left(B_{\gamma t/\gamma}^H - \int_0^{\gamma t} e^{-(\gamma t - s)} B_{s/\gamma}^H\,ds\right)_{t\geq 0} \stackrel{d}{=} \gamma^{-H}(B_{1,\gamma t}^H)_{t\geq 0}.$$

Next, we represent the vector $\tau_H(\gamma)$ in terms of the random quantities on the left-hand side in the last display. For $\gamma > 1$, set $s_H(\gamma) = \gamma g_H(\gamma)$ and

$$V_\gamma^H = \gamma^{2H-1}\left(\tfrac{1}{2}(B_{\gamma,1}^H)^2 + \gamma \int_0^1 (B_{\gamma,s}^H)^2\,ds\right), \qquad W_\gamma^H = \gamma^{2H}\int_0^1 (B_{\gamma,s}^H)^2\,ds.$$

In view of (2.11) and (4.21), we have

$$\operatorname{diag}(s_H(\gamma), s_H(\gamma), 1)(V_\gamma^H - \mu_H^{-1}, W_\gamma^H - \mu_H^{-1}, \tau_{H,4}(\gamma) - \gamma)'$$
$$(4.22) \quad \stackrel{d}{=} g_H(\gamma)\int_0^\gamma [(B_{1,s}^H)^2 - \mu_H^{-1}]\,ds\begin{pmatrix}1\\1\\0\end{pmatrix} + \frac{1}{2}\frac{(B_{1,\gamma}^H)^2}{\gamma^{-1}\int_0^\gamma (B_{1,s}^H)^2\,ds}\begin{pmatrix}0\\0\\1\end{pmatrix}$$
$$+ \frac{1}{2}g_H(\gamma)(B_{1,\gamma}^H)^2\begin{pmatrix}1\\0\\0\end{pmatrix}, \qquad \gamma > 1.$$

Note that $g_H(\gamma) \to 0$ and $\gamma^{-1}\int_0^\gamma (B_{1,u}^H)^2\,du \stackrel{P}{\to} \mu_H^{-1}$, as $\gamma \to \infty$, by Lemma 4.2. Consequently, we get from Lemma 4.2 that, as $\gamma \to \infty$, the left-hand side in (4.22) converges jointly in distribution to the random vector

$$(4.23) \qquad (\sigma_H X^H, \sigma_H X^H, Y^2/2)'.$$

Applying the delta method to $\psi(v,w) = (w^{-1/2}, w^{-1}, vw^{-1/2})'$, as $\gamma \to \infty$,

$$s_H(\gamma)[\mathrm{diag}(\gamma^{-H}, \gamma^{-2H}, \gamma^{H-1})(\tau_{H,1}(\gamma), \tau_{H,2}(\gamma), \tau_{H,3}(\gamma))' - \psi(\mu_H^{-1}, \mu_H^{-1})]$$

$$= s_H(\gamma)[\psi(V_\gamma^H, W_\gamma^H) - \psi(\mu_H^{-1}, \mu_H^{-1})] \xrightarrow{d} \sigma_H X^H \psi'(\mu_H^{-1}, \mu_H^{-1}) \begin{pmatrix} 1 \\ 1 \end{pmatrix}.$$

As we have joint convergence in (4.22) to the vector in (4.23), the proof of Theorem 3.1 is completed by elementary calculus.

4.4. *Proof of Theorem 3.2.* Fix $H \in (0,1)$. For all $\gamma < 0$, we set

(4.24)
$$T_{H,1}(\gamma) = \tfrac{1}{2}(B_{-1,|\gamma|}^H)^2 - \int_0^{|\gamma|} (B_{-1,u}^H)^2 \, du,$$

$$T_{H,2}(\gamma) = \int_0^{|\gamma|} (B_{-1,u}^H)^2 \, du.$$

It follows from (2.5) that the following identity holds for the finite dimensional distributions:

$$(B_{-\gamma,t}^H)_{t\geq 0} \stackrel{d}{=} |\gamma|^{-H} (B_{-1,|\gamma|t}^H)_{t\geq 0}, \qquad \gamma < 0, \; H \in (0,1).$$

Consequently, for $\gamma < 0$,

$$\mathrm{diag}(|\gamma|^{-(2H+1)/2} e^{|\gamma|}, |\gamma|^{-2H-1} e^{2|\gamma|}, |\gamma|^{(2H-1)/2}, |\gamma|^{-1} e^{|\gamma|}) \tau_H(\gamma)$$

$$\stackrel{d}{=} \mathrm{diag}(e^{|\gamma|}, e^{2|\gamma|}, 1, e^{|\gamma|})$$

$$\times \left( \frac{1}{T_{H,2}(\gamma)^{1/2}}, \frac{1}{T_{H,2}(\gamma)}, \frac{T_{H,1}(\gamma)}{(T_{H,2}(\gamma))^{1/2}}, \frac{T_{H,1}(\gamma)}{T_{H,2}(\gamma)} \right)'.$$

Integrating by parts and summing and subtracting terms, for $\gamma < 0$, we have

(4.25)
$$e^\gamma T_{H,1}(\gamma) = \left( \int_0^\infty e^{-s} B_s^H \, ds \right) \left( \int_0^{|\gamma|} e^{\gamma+s}(B_{|\gamma|}^H - B_s^H) \, ds \right)$$
$$+ e^\gamma \tilde{T}_{H,1}(\gamma),$$

where we set

$$\tilde{T}_{H,1}(\gamma) = \int_0^{|\gamma|} e^s B_s^H \int_s^\infty e^{-t} B_t^H \, dt \, ds + \tfrac{1}{2}(B_{|\gamma|}^H)^2 - \int_0^{|\gamma|} (B_s^H)^2 \, ds$$
$$- e^{|\gamma|} B_{|\gamma|}^H \int_{|\gamma|}^\infty e^{-s} B_s^H \, ds + B_{|\gamma|}^H \int_0^\infty e^{-s} B_s^H \, ds, \qquad \gamma < 0.$$

As $B_t^H/t \to 0$ for $t \to \infty$ and $H \in (0,1)$, almost surely, we note that $C_H = \sup_{t\geq 0} \frac{|B^H(t)|}{1+t}$ is finite, almost surely. Consequently, we have the following

chain of inequalities [$\gamma < 0$], almost surely,

$$\left|\int_0^{|\gamma|} e^s B_s^H \int_s^\infty e^{-t} B_t^H \, dt \, ds\right| \leq C_H^2 \int_0^{|\gamma|} e^s(1+s) \int_s^\infty e^{-t}(1+t) \, dt \, ds$$

$$= C_H^2 \int_0^{|\gamma|}(1+s) \int_0^\infty e^{-t}(1+s+t) \, dt \, ds \leq C' C_H^2 (1+|\gamma|^3),$$

where $C' > 0$ is a finite deterministic constant not depending on $\gamma$.

Exploring similar bounds for the other terms in (4.26) yields, almost surely,

$$|\tilde{T}_{H,1}(\gamma)| \leq C'' \max(C_H, C_H^2)(1+|\gamma|^3),$$

where $C''$ is a finite and deterministic constant not depending on $\gamma$. Consequently, almost surely,

(4.26) $$e^\gamma \tilde{T}_{H,2}(\gamma) \to 0, \qquad \gamma \to -\infty.$$

Similar,

(4.27) $$\tfrac{1}{2} e^{2|\gamma|} \left(\int_0^\infty e^{-s} B_s^H \, ds\right)^2 + \tilde{T}_{H,3}(\gamma), \qquad \gamma < 0,$$

where $e^{2\gamma} \tilde{T}_{H,3}(\gamma) \to 0$, for $\gamma \to -\infty$, almost surely.

Finally, for all $H \in (0,1)$ and $\gamma \to -\infty$, note that

(4.28) $$\begin{pmatrix} \int_0^{|\gamma|} e^{\gamma+s}(B_{|\gamma|}^H - B_s^H) \, ds \\ \int_0^\infty e^{-s} B_s^H \, ds \end{pmatrix} \xrightarrow{d} \left(\frac{\Gamma(2H+1)}{2}\right)^{1/2} \begin{pmatrix} Y \\ Z \end{pmatrix},$$

where $(Y, Z)$ is standard normal. Indeed, as all corresponding vectors are centered Gaussian, it suffices to investigate their covariance matrices. Clearly, $\int_0^\infty e^{-s} B_s^H \, ds$ is Gaussian with a mean zero and variance $\Gamma(2H+1)/2$. Furthermore,

$$\int_0^{|\gamma|} e^{\gamma+s}(B_{|\gamma|}^H - B_s^H) \, ds \stackrel{d}{=} \int_0^{|\gamma|} e^{-s} B_s^H \, ds, \qquad \gamma < 0.$$

The right-hand side tends to $\int_0^\infty e^{-s} B_s^H \, ds$, almost surely and in mean square as $\gamma \to -\infty$.

For all $\gamma \in \mathbb{R}$ and $H \in (0,1)$, we have the following identity,

$$E\left(\int_0^\infty e^{-t} B_t^H \, dt\right) \left(\int_0^{|\gamma|} e^{\gamma+s}(B_{|\gamma|}^H - B_s^H) \, ds\right)$$

$$= \int_0^\infty e^{-t} \int_0^{|\gamma|} e^{-s} E B_t^H (B_{|\gamma|}^H - B_{|\gamma|-s}^H) \, ds \, dt.$$

In view of (2.4), expanding the integrand into powers of $s/|\gamma|$, $t/|\gamma|$ and $(s+t)/|\gamma|$ up to the second order, this shows that the integrand tends to zero for $\gamma \to -\infty$ and $s,t \in \mathbb{R}$. Note that $|EB_t^H(B_{|\gamma|}^H - B_{|\gamma|-s}^H)| \leq |t|^{2H}|s|^{2H}$ by means of the Cauchy–Schwarz inequality. Thus, the dominated convergence theorem is applicable to the right-hand side in the last display giving (4.28) for all $H \in (0,1)$ and $\gamma \to -\infty$. In view of (4.25), (4.26) and (4.27), (i) now follows from (4.28).

## APPENDIX: THE DIAGRAM FORMULA

We refer to Bollobás [5] for basic notation and facts of graph theory. The sets of vertices and edges of an undirected graph $D$ are denoted by $\mathcal{V}(D)$ and $\mathcal{E}(D)$, respectively. Throughout, we identify $\mathcal{E}(D)$ with a subset of the set $\{\{v,w\} : v,w \in \mathcal{V}(D), v \neq w\}$. The *degree* of a vertex $v \in \mathcal{V}(D)$ is defined to be the cardinality of the set $\{e \in \mathcal{E}(D) : v \in e\}$.

In the sequel, let $p, q \in \mathbb{N}$. A *diagram $D$ with $p$ levels and rows of length $q$* is an undirected graph with vertex set $\mathcal{V}(D) = \{1, \ldots, p\} \times \{1, \ldots, q\}$, satisfying the following properties:

(i) If $\{(l_1, k_1), (l_2, k_2)\} \in \mathcal{E}(D)$ then $l_1 \neq l_2$.
(ii) Each vertex has degree one.

The (possibly empty) set of all diagrams with $p$ levels and rows of length $q$ is denoted by $\mathcal{D}(p, q)$. For $D \in \mathcal{D}(p, q)$ and $e = \{(k, m), (l, n)\} \in \mathcal{E}(D)$, set

$$(A.1) \quad m(e) = m_D(e) = \min\{k, l\}, \qquad M(e) = M_D(e) = \max\{k, l\}.$$

Now let $\mathcal{H}_q$ be the Hermite polynomial of order $q$ in (1.2) and $(X_1, \ldots, X_p)$ be a Gaussian vector with $EX_j = 0$ and $EX_j^2 = 1$ ($1 \leq j \leq p$). It follows from the *diagram formula* (cf. Breuer and Major [7]) that we have

$$(A.2) \quad E\prod_{j=1}^p \mathcal{H}_q(X_j) = \sum_{D \in \mathcal{D}(p,q)} \prod_{e \in \mathcal{E}(D)} E[X_{m(e)} X_{M(e)}],$$

with the conventions $\prod_\varnothing = 1$ and $\sum_\varnothing = 0$.

**Acknowledgments.** We would like to thank the Editor, an Associate Editor and two anonymous referees for their insightful suggestions and relevant references. We would also like to thank C. C. Heyde and R. A. Maller for many helpful discussions during the project. In particular, we would like to dedicate this paper to Professor Chris Heyde, who recently passed away due to a long-term illness. Part of this research was completed when the first author visited the Institute of Mathematical Sciences (IMS) at CUHK in 2004 and 2006. Research support from the IMS is gratefully acknowledged.

School of Mathematical Sciences  
Monash University  
VIC 3800  
Australia  
E-mail: Buchmann@sci.monash.edu.au

Department of Statistics  
Chinese University of Hong Kong  
Shatin, New Territories  
Hong Kong  
E-mail: nhchan@sta.cuhk.edu.hk